\documentclass[a4paper,10pt]{amsart}
\usepackage{textcomp}
\usepackage{amssymb, amsmath, amsthm}
\usepackage{xypic}
\usepackage{enumerate}

\xyoption{all}

\theoremstyle{theorem}
\newtheorem{theorem}{Theorem}[section]
\newtheorem{lemma}[theorem]{Lemma}
\newtheorem{proposition}[theorem]{Proposition}

\theoremstyle{definition}
\newtheorem{definition}[theorem]{Definition}
\newtheorem{example}[theorem]{Example}

\theoremstyle{remark}
\newtheorem{remark}[theorem]{Remark}

\def\m{{\mathfrak {m}}}
\def\A{{\mathcal {A}}}
\def\C{{\mathcal {C}}} 
\def\D{{\mathcal {D}}}
\def\triv{\stackrel{\sim}{\ra}}
\def\ra{\rightarrow}
\def\R{R}
\def\Z{{\mathbb  {Z}}}

\title{Classification of certain cellular classes of chain complexes}

\begin{document}

\author[Kiessling] 
{Jonas Kiessling} 
 
\address{Kiessling: Department of Mathematics\\ 
The Royal Institute of Technology\\ 
S -- 100 44 Stockholm\\ 
SWEDEN} 
 
\email{jonkie@kth.se} 
 
\thanks{Research supported by grant KAW 2005.0098 from 
the Knut and Alice Wallenberg Foundation.} 

\begin{abstract}
Let $(R, \m)$ be a local commutative ring. Suppose that $\m$ is principal and that $\m^2 = 0$. We give a complete description of the cellular lattice of perfect chain complexes of modules over this ring.  
\end{abstract}

\maketitle 

\section{Introduction}
An explicit classification of thick subcategories of finite spectra (~\cite{MR1652975} and ~\cite{MR1333942})  has been an important achievement of homotopy theory. This stable classification has been used to give an explicit classification of unstable Bousfield classes of finite suspension spaces~\cite{MR1397720}.  In contrast, an analogous classification of cellular classes of finite (suspension) spaces is out of reach as that would lead  to a classification of ideals in the stable homotopy groups of spheres. 

Thick subcategories of compact objects in the derived category of a ring are also well understood~\cite{MR1174255}. Recently Bousfield classes (or Acyclic classes) of chain complexes have been classified, see \cite{Staney} or \cite{Kiessling}. However, as for spaces, the classification of cellular classes is more subtle. The aim of this paper is to give some examples of rings for which a classification can be obtained. For some properties of cellular classes of chain complexes over a Noetherian ring, see \cite{Kiessling}. 

Throughout this paper all chain complexes are non-negatively graded chain complexes of modules over some fixed commutative ring.

\begin{definition} \label{def:cell1}
Fix a chain complex $A$.
We let $\C(A)$ denote the smallest collection of chain complexes satisfying the following properties: 1. The collection $\C(A)$ contains $A$, 2. It is closed under arbitrary sums and weak equivalences (i.e. homology isomorphisms) 3. If $0 \ra X \ra Y \ra Z \ra 0$ is an exact sequence with $X,Y \in \C(A)$ then also $Z \in \C(A)$. If $X \in \C(A)$ then we write that $X \gg A$ and say that $X$ is $A-$cellular.
\end{definition}

In the paper we describe explicitly the cellular relations between perfect chain complexes of modules
over a local ring $(R,m)$  such that  $m$ is principal and $m^2=0$. 

For such a ring and numbers $i, j\geq 0$, we let:
\[
(\Sigma^i E_j)_{n}= \left\{ \begin{array}{cl}
R & \textrm{for $i \leq n \leq i+j$} \\
0 & \textrm{otherwise} \end{array} \right.
\]
with the differentials given by multiplication by $(-1)^i r$, where $r$ is some generator of the maximal ideal $\m$.  The isomorphism type of $\Sigma^i E_j$ does not depend on the choice of the generator $r$.

We prove that:

\begin{theorem}
Let  $(R,m)$  be a local ring such that  $m$ is principal and $m^2=0$.
Let $A$ be a perfect chain complex that is not weakly equivalent to 0. Then there exists $(i,j)$ such that:
\[
\C(A) = \C(\Sigma^i E_j)
\]
Moreover $\Sigma^{j'} E_{i'} \gg \Sigma^{j} E_{i}$ if and only if either $ i' \geq i$ or $i' = i$ and $j' \geq j$.
\end{theorem}

The key element in the proof of the theorem is a classification of the perfect chain complexes.
We show that any perfect chain complex $A$ splits into a sum of a contractible chain complex and a sum of $\Sigma^i E_j$'s (see Lemma \ref{lemma:splitting}).

\section{Notation} \label{section:notation}

We let $R$ denote some arbitrary commutative ring. By a \emph{chain complex} $X$ we mean a non-negatively graded chain complex of $R$-modules. We use the homological grading, i.e. the differential of $X$ lowers the degree. Recall that the category of chain complexes of $R$-modules $Ch_{\geq 0}(R)$ is a model category~\cite{MR1361887}. In this model category 
a \emph{weak equivalence} is a map that induces an isomorphism on homology. 
A \emph{cofibration} is an injective map such that the cokernel is projective in each degree. A \emph{fibration} is map which is surjective in all positive degrees. We let $\triv$ denote a weak equivalence. A cofibrant chain complex $X$ is a chain complex such that the canonical map $0 \ra X$ is a cofibration, or explicitly, it is a chain complex of projective modules.

If $f:X \ra Y$ is any map of chain complexes then we can factor $f = f'' f'$, where $f':X \ra X'$ is a cofibration and $f'':X' \triv Y$ is a weak equivalence (\cite{MR1361887}). Hence, any map is a cofibration up to a weak equivalence. 

A complex $X$ is called \emph{perfect} if it is cofibrant and $\oplus_i X_i$ is finitely generated. 

We let $Hom$ denote the \emph{hom-complex}. It is defined as follows: if $X, Y \in Ch_{\geq 0}(\R)$ then $Hom(X,Y)_n = \prod_i hom(X_i, Y_{i+n})$ for $n>0$ and $Hom(X,Y)_0$ is the set of maps of chain complexes from $X$ to $Y$ with the induced $R$-module structure. The differential takes $\{f_i:X_i \ra Y_{i+n}\}_i \in Hom(X,Y)_n$ to $\{\partial f_i + (-1)^n f_{i-1} \partial \}_i$. If $A$ is cofibrant, then $Hom(A, \bullet)$ preserves weak equivalences and fibrations (a consequence of Brown's lemma, see \cite{MR1361887}).  

We let $\Sigma^i$ denote the \emph{shift operator} i.e $(\Sigma^i X)_j = X_{j-i}$ and $\partial_{\Sigma^i} = (-1)^i \partial$. The \emph{cone} of a map $f:X \ra Y$ is a chain complex $C(f)$ defined by; $C(f)_n = Y_n \oplus X_{n-1}$. The differential of $C(f)$ maps $(y,x) \in C(f)_n$ to $(\partial^Y(y) + f(x), -\partial^X(x))$. Note that there is a canonical map $Y \ra C(f)$ and the cokernel of this map is isomorphic to $\Sigma^1 X$. If $0 \ra X \ra Y \ra Z \ra 0$ is an exact sequence of chain complexes then there is a natural map from the cone of $X \ra Y$ to $Z$. From the induced long exact sequences of homology, this map is a weak equivalence.

A more detailed account of the theory of chain complexes can be found for instance in \cite{MR1269324}.

Associated to any module $M$ are the \emph{sphere complex}, $S^n(M)$, and the \emph{disk complex}, $D^n(M)$, defined by:
\[
(S^n(M))_i = \left\{ \begin{array}{cl}
M & i = n \\
0 & i \neq n \end{array} \right.
\quad 
(D^n(M))_i = \left\{ \begin{array}{cl}
M & \textrm{$i = n$ or $i = n-1$} \\
0 & \textrm{otherwise} \end{array} \right.
\]
With the differential $\partial_n = 1_{R}$ in $D^n(M)$. For short we let $S^n := S^n(R)$ and $D^n := D^n(R)$.

\section{Cellular relation}
Recall that there is an alternative description of cellularity via a universal property (see \cite{MR1392221}):

\begin{proposition} \label{prop:cell2}
Let $X$ and $A$ be cofibrant chain complexes. Then $X$ is $A$-cellular if and only if for all maps $f$ such that $Hom(A,f)$ is a weak equivalence, the map $Hom(X,f)$ is also a weak equivalence.
\end{proposition}

The cellular relation is transitive, in other words if $X \gg A$ and $Y \gg X$ then $Y \gg A$. 

To determine whether a given chain complex $X$ belongs to $\C(A)$ is in general a hard question. In Proposition~\ref{prop:surjective} we give a workable criteria for detecting cellularity for a very particular choice of ring $R$. We now list certain properties of cellularity.

\begin{proposition} \label{prop:properties}
\begin{enumerate}[(i)]
\item All chain complexes are $S^0$-cellular.
\item If $X$ is acyclic (i.e. all homology group vanish), then $X \gg A$ for any chain complex $A$.
\item The collection $\C(A)$ is closed under retracts.
\item Suppose that $A$ is a cofibrant chain complex. If $X\gg A$ and $H_0(A) \neq 0$ then there is a set $I$ and a map $f: \oplus_I A \ra X$ such that $H_0(f)$ is surjective.
\item If $X \gg A$ then $\Sigma^n X \gg A$ for all $n \geq 0$.
\item $\Sigma^1 X \gg \Sigma^1 A \Leftrightarrow X \gg A$.
\end{enumerate}
\end{proposition}

\begin{proof}
We shall only give an outline of the proof. Statement $(i)$ follows from the isomorphism $Hom(S^0, X) \cong X$ and \ref{prop:cell2}. To prove the second statement note first that $0 \in \C(A)$ since $Hom(0,Y) = 0$. If $X$ is acyclic then $0 \ra X$ is a weak equivalence and since $\C(A)$ is closed under weak equivalences $X \gg A$. A retract of an isomorphism is an isomorphism so $(iii)$ is a consequence of \ref{prop:cell2}. 

We fix a cofibrant chain complex $A$ and let $\D$ denote the collection of all chain complexes $X$ such that there is a set $I$ and a map $f: \oplus_I A \ra X$ surjective on $H_0$. It is a standard result in homological algebra (\cite{MR1269324}, p. 388) that if $f:X \triv Y$ is a weak equivalence and $g:A \ra Y$ is any map, then because $A$ is cofibrant, there is a map $h:A \ra X$ such that $g$ and $fh$ are homotopic. As a consequence, the collection $\D$ is closed under weak equivalences. It is also closed under sums. Finally if $0 \ra X \ra Y \ra Z \ra 0$ is an exact sequence such that $X,Y \in \D$ then also $Z \in \D$. By definition $\C(A)$ was the smallest collection satisfying these properties, hence $\C(A) \subset \D$ and we have proved $(iv)$.

Statements $(v)$ and $(vi)$ are direct consequences of \ref{prop:cell2}. 

\end{proof}

\section{Two out of three property}

We say that a collection of chain complexes $\C$ satisfies the \emph{two out of three property} if given any exact sequence $0 \ra X \ra Y \ra Z \ra 0$ such that two out of $X$, $Y$ and $Z$ belong to $\C$ then so does the third. The collection $\C(A)$ does not in general satisfy the two out of three property. For instance it follows from \ref{prop:properties} that $\C(S^1)$ equals the collection of all chain complexes $X$ such that $H_0 X = 0$. There is an exact sequence $0 \ra S^0 \ra D^1 \ra S^1 \ra 0$ and $D^1$ is $S^1$-cellular, but $S^0$ is not.

Collections of possibly unbounded chain complexes of modules over a Noetherian ring, satisfying the two out of three property, that are closed under sums and weak equivalences have been classified by Neeman in \cite{MR1174255}. They are in 1-1 correspondence with arbitrary sets of prime ideals in $R$. 

A collection $\C$ of chain complexes is closed under extensions if given any exact sequence $0 \ra X \ra Y \ra Z \ra 0$ with $X, Z \in \C$ then also $Y \in \C$. The collection $\C(A)$ is in general not closed under extensions. In analogy with cellularity we can now define a relation called \emph{acyclicity}:

\begin{definition}
Fix a chain complex $A$. Let $\A(A)$ denote the smallest collection of chain complexes satisfying the following properties: 1. The collection $\A(A)$ contains $A$, 2. It is closed under arbitrary sums, 3. If in an exact sequence $0 \ra X \ra Y \ra Z \ra 0$, either $X$ and $Y$ or $X$ and $Z$ belong to $\A(A)$, then so does the third. If $X \in \A(A)$ then we write that  $X > A$ and say that $X$ is $A$-\emph{acyclic}.
\end{definition}

The collection $\A(A)$ is in particular closed under extensions. By definition $C(A) \subset \A(A)$. This inclusion is in general strict:
\begin{example} 
Recall from the introduction the definition of $E_i := \Sigma^0 E_i$. Up to a weak equivalence, $E_1$ is an extension of $E_3$ by $E_2$: Let $r$ denote some generator of the maximal ideal. Multiplication by $r$ in degree 0 and the zero map in higher degrees defines a map $f:E_2 \ra E_1$. The cone of $f$ is isomorphic to $E_3$. This gives the sequence: $E_2 \ra  E_1 \ra E_3$

We later show (Theorem \ref{theorem:cellularity}) that $E_3 \gg E_2$ (that is $E_3$ is $E_2$-cellular) and that $E_1$ is \emph{not} $E_2$-cellular. However $E_1 > E_2$ since $E_1$ is an extension of $E_3$ by $E_2$ and $E_3 \gg E_2$.
\end{example}

The relation $>$ between perfect chain complexes of modules over a Noetherian ring is well understood. In fact $X>A$ if and only if for every $n$ the following holds: if $p \subset R$ is a prime ideal such that $X_n \otimes R_p \neq 0$ then $A_i \otimes R_p \neq 0$ for some $i \leq n$. The proof of this result can be found in the papers \cite{Staney} and \cite{Kiessling}. 

The following Proposition establishes an important connection between cellularity and acyclicity. An analogous result for topological spaces was obtained by Dror-Farjoun in \cite{MR1392221}. 

\begin{proposition} \label{prop:connection}
Fix a short exact sequence $0 \ra X \ra Y \ra Z \ra 0$. If $X \gg A$ and $Z > \Sigma^1 A$ then $Y \gg A$.
\end{proposition}
\begin{proof}
See \cite{Kiessling}. A proof of the topological analog can be found in \cite{MR1392221}.  
\end{proof}

\section{A Key Lemma}

For the rest of this paper we fix a local commutative ring $R$ with a principal maximal ideal $\m$ such that $\m^2 = 0$. We let $k$ denote the residue field $R / \m$. We choose some generator $r$ of $\m$. Note that in such a ring all non-unitary elements are of the form $x = r' r$, for some unit $r'$. Good examples to keep in mind are $R = \Z / (p^2)$ ($p$ prime) and $R = k[X]/(X^2)$, for some field $k$. 

\begin{definition} An injective map $D^n \ra X$ is called an \emph{embedded disk}. \end{definition}

Recall from the introduction the special class of chain complexes:
\[
(E_j)_n = \left\{ \begin{array}{cl}
R & \textrm{if $n \leq j$} \\
0 & \textrm{if $n > j$} \end{array} \right.
\]
With $\partial_n: (E_j)_n \ra  (E_j)_{n-1}$ multiplication by $r$.
For instance, $E_{\infty}$ is a projective resolution of $k$ and $E_1 = S^0$.

\begin{lemma} \label{lemma:splitting}
Let $X$ be any perfect chain complex. Then there is a splitting of $X$:
\[
X \cong P \oplus Q
\]
where $P$ is acyclic and $Q$ can be written as a finite sum
\[
Q \cong \oplus \Sigma^i E_j
\]
\end{lemma}
\begin{proof} 
The proof is divided into several steps and takes up the rest of this section. Fix some perfect chain complex $X$.

\emph{Step 1}

We first split off the contractible part of $X$.
\begin{remark} 
$R$ is injective as a module over itself. Hence $D^n$ is an injective object, i.e. any embedded disk $D^n \ra X$ is split. 
\end{remark}
A consequence of this remark is that $X$ will split into a direct sum
\[
X \cong P \oplus \tilde{X}
\]
where $P$ is an acyclic complex (a finite sum of disks) and $\tilde{X}$ has no embedded disks.

\emph{Step 2}

We can now assume that $X$ has no embedded disks. We also want to assume that $X_0 \neq 0$. 

\begin{remark} \label{remark:disk}
For a perfect complex $X$, to have no embedded disks is equivalent to $\partial_n(X_n) \subset r X_{n-1}$ for all $n$. Such a complexes are also known as \emph{minimal}.
\end{remark}

Suppose that $H_0(X) = 0$. Then $\partial_1$ is surjective. From the remark we see that $X_0 = 0$. Hence $X = \Sigma^1 Y$ for some perfect complex $Y$.

\emph{Step 3}

By step 1 and 2 we can assume that $X$ is a perfect chain complex, containing no embedded disks and that $X_0 \neq 0$. 

Since $H_0(X) \neq 0$ there is a surjection $X \ra S^0(k)$ and because $X$ is cofibrant this map factors:
\[
X \ra E_{\infty} \triv S^0(k)
\]
where $E_{\infty} \triv S^0(k)$ is an acyclic fibration (i.e. a fibration and a weak equivalence).

The chain complex $X$ is perfect, in particular $\oplus_i X_i$ is finitely presented, so the map $X \ra E_{\infty}$ factors through some $E_n$. We note that the map $X \ra E_n$ is surjective in degree 0. There is some smallest $n_0$ such that there is a map $f:X \ra E_{n_0}$ surjective in degree 0. We claim that $f$ is surjective and has a section, so that $X \cong E_{n_0} \oplus \tilde{X}$.

First we show that $f$ is surjective. Suppose that this is not the case. Then there is some $m \leq n_0$ such that im $f_m \subset $ ker $\partial_m$. We truncate $f$ at $m-1$ and get a map $\tilde{f}:X \ra E_{m-1}$ with $\tilde{f}_j = f_j$ for $j \leq m-1$ and $\tilde{f}_j = 0$ for $j \geq m$. Then $\tilde{f}$ is surjective in degree 0, contradicting the minimality of $n_0$. Hence $f$ is surjective.

We fix a generator $e_j$ of $(E_{n_0})_j$ for each $j$. From the surjectivity of $f_{n_0}$ it follows that there is some $x_{n_0} \in X_{n_0}$ such that $f_{n_0}(x) = e_{n_0}$. By remark~\ref{remark:disk} there is some $x_{n_0 - 1}$ such that $r x_{n_0 - 1} = \partial(x_{n_0})$. Moreover $f_{n_0 - 1}(x_{n_0 - 1}) = a_{n_0 - 1} e_{n_0 - 1}$, for some unit $a_{n_0 - 1}$. Inductively we obtain a sequence of elements $(x_0, \ldots , x_{n_0})$ and $(a_0, \ldots, a_{n_0})$. We define a map $s: E_{n_0} \ra X$ by $s_j(e_j) = x_j$, $j \leq n_0$. This is well defined since $\partial \circ s_j (e_j) = r x_{j-1} = s_{j-1}(\partial (e_{j-1}))$. By construction $(f \circ s)_j (e_j) = a_j e_j$ so that $f \circ s$ is an isomorphism.

This determines a splitting of $X$ into
\[
X \cong E_{n_0} \oplus \tilde{X}
\]

We can repeat the above discussion with $\tilde{X}$ instead of $X$. This yields the required splitting formula.

\end{proof}

\section{Statement of Results} \label{section:main}

Recall that $R$ is assumed to be a commutative local ring with a principal maximal ideal $\m$ such that $\m^2 = 0$. We are now in a position to give a complete description of the cellular lattice of perfect complexes.

First we look at  the cellularity relations among  $E_i$'s. We begin with an observation.

\begin{remark}\label{remark:map}
If $i > j$ then there is no map $f: E_i \ra E_j$, such that $H_0 (f)$ is non-zero.
\end{remark}
This remark in combination with the following proposition is enough to classify the cellular relations between the $E_i$'s.

\begin{proposition}\label{prop:surjective}
Let $X$ be a perfect complex such that $H_0(X) \neq 0$. Then $Y$ is $X-$cellular if and only if there is a set $I$ and a map $f: \oplus_{i \in I} X \ra Y$ such that $f$ induces an epimorphism on the $H_0$.
\end{proposition}
\begin{proof}
We first fix some perfect complex $X$ with $H_0(X) \neq 0$. Let $D(X)$ denote the class of all complexes $Y$ such that there is some $I$ and $f:\oplus_I Y \ra X$ with $H_0(f)$ surjective. The statement is that $D(X) = C(X)$.

By Proposition~\ref{prop:properties}. To prove that $D(X) \subset C(X)$ we first note that since $X$ is cofibrant, $X \otimes \bullet$ preserves cellularity: if $A \gg B$ then $A \otimes X \gg B \otimes X$. Here $\otimes$ denotes the ordinary tensor product of chain complexes (see \cite{MR1269324}). Since $X \otimes S^0 \cong X$ and $S^0(k) \gg S^0$ we can conclude that $X \otimes S^0(k) \gg X$. We can always assume that there are no embedded disks in $X$ and in this case $H_0(X \otimes k) \cong (X \otimes k)_0$. By assumption $H_0(X) \neq 0$, so there is a retraction: $S^0(k) \ra X \otimes k \ra S^0(k)$. It follows from Proposition~\ref{prop:properties} that $S^0(k) \gg X \otimes k \gg X$. 

Fix a map $f: \oplus_{i \in I} X \ra Y$ such that $H_0(f)$ is surjective. We can assume that $f: \oplus_{i \in I} X \ra Y$ is a cofibration (see section \ref{section:notation}). Let $Z := Y / f(\oplus_{i \in I} X)$. We wish to show that $Y$ is $X-$cellular. By Proposition~\ref{prop:connection} it is enough to show that $Z > \Sigma^1 X$.  

We assumed that $H_0(f_0)$ is surjective, so $H_0(Z) \cong 0$. As a consequence $Z \gg S^1$. There is an isomorphism of $R-$modules $\m \cong k$. The exact sequence $0 \ra k \ra R \ra k \ra 0$ shows that $S^0 > S^0(k)$ or equivalently that $S^1 > S^1(k)$. In the paragraph above we showed that $S^1(k) \gg \Sigma^1 X$. In all we have that $Z \gg S^1 > S^1(k) \gg \Sigma^1 X$, i.e. $Z > \Sigma^1 X$. This concludes the proof of the proposition. 

\end{proof}    

Consider the set $S$ of all pairs $(i,j)$ of non-negative integers. We order this set by declaring $(i,j) < (i', j')$ if either $i < i'$ or $i = i' $ and $j < j'$. To each perfect complex $X$ we assign the subset $S_X$ of $S$ consisting of all pairs $(i,j)$ such that $\Sigma^i E_j$ appears in the splitting of $X$ as in lemma~\ref{lemma:splitting}. Finally we let $(i_X,j_X)$ denote the minimal element of $S_X$. We have shown that:

\begin{theorem}\label{theorem:cellularity}
Let  $(R,m)$  be a local ring such that  $m$ is principal and $m^2=0$. Let $X$ be a perfect complex over $R$. Then 
\[
\Sigma^{i'} E_{j'} \gg \Sigma^i E_j
\]
if and only if $(i,j) \leq (i',j')$. Moreover: 
\[
\C(X) = \C(\Sigma^{i_X} E_{j_X})
\]
unless $X$ is contractible, in which case $\C(X) = \C(0)$.
\end{theorem}




\begin{thebibliography}{Bou96}

\bibitem[Bou96]{MR1397720}
A.~K. Bousfield, \emph{Unstable localization and periodicity}, Algebraic
  topology: new trends in localization and periodicity (Sant Feliu de Gu\'\i
  xols, 1994), Progr. Math., vol. 136, Birkh\"auser, Basel, 1996, pp.~33--50.

\bibitem[DH95]{MR1333942}
Ethan~S. Devinatz and Michael~J. Hopkins, \emph{The action of the {M}orava
  stabilizer group on the {L}ubin-{T}ate moduli space of lifts}, Amer. J. Math.
  \textbf{117} (1995), no.~3, 669--710.

\bibitem[DS95]{MR1361887}
W.~G. Dwyer and J.~Spali{\'n}ski, \emph{Homotopy theories and model
  categories}, Handbook of algebraic topology, North-Holland, Amsterdam, 1995,
  pp.~73--126.

\bibitem[Far96]{MR1392221}
Emmanuel~Dror Farjoun, \emph{Cellular spaces, null spaces and homotopy
  localization}, Lecture Notes in Mathematics, vol. 1622, Springer-Verlag,
  Berlin, 1996.

\bibitem[HS98]{MR1652975}
Michael~J. Hopkins and Jeffrey~H. Smith, \emph{Nilpotence and stable homotopy
  theory. {II}}, Ann. of Math. (2) \textbf{148} (1998), no.~1, 1--49.

\bibitem[Kie]{Kiessling}
Jonas Kiessling, \emph{Properties of cellular classes}, In preparation.

\bibitem[Nee92]{MR1174255}
Amnon Neeman, \emph{The chromatic tower for {$D(R)$}}, Topology \textbf{31}
  (1992), no.~3, 519--532, With an appendix by Marcel B\"okstedt.

\bibitem[Sta]{Staney}
Don Stanley, \emph{Invariants of t-structures and classification of nullity
  classes}, arXiv:math/0602252v1.

\bibitem[Wei94]{MR1269324}
Charles~A. Weibel, \emph{An introduction to homological algebra}, Cambridge
  Studies in Advanced Mathematics, vol.~38, Cambridge University Press,
  Cambridge, 1994.

\end{thebibliography}

\providecommand{\bysame}{\leavevmode\hbox to3em{\hrulefill}\thinspace}
\providecommand{\MR}{\relax\ifhmode\unskip\space\fi MR }
\providecommand{\MRhref}[2]{%
  \href{http://www.ams.org/mathscinet-getitem?mr=#1}{#2}
}
\providecommand{\href}[2]{#2}

\end{document}